\documentclass[a4paper,12pt]{amsart}
\usepackage{amssymb}
\author{D.~A.~Timashev}
\thanks{Supported by the NATO research scholarship 350590A}
\title[Complexity and growth of multiplicities]
{Complexity of homogeneous spaces\\
and growth of multiplicities}
\date{May 29, 2003}
\address{Department of Higher Algebra\\
Faculty of Mechanics and Mathematics\\
Moscow State University\\
119992 Moscow\\
Russia}
\curraddr{Institut Fourier\\
UFR des Math{\'e}matiques\\
UMR 5582\\
B.~P.~74\\
38402 Saint-Martin d'H{\`e}res CEDEX\\
France}
\email{timashev@mech.math.msu.su}
\urladdr{http://mech.math.msu.su/department/algebra/staff/timashev}

\newcommand{\kk}{\Bbbk}
\newcommand{\PP}{\mathbb{P}}
\newcommand{\ZZ}{\mathbb{Z}}
\newcommand{\codim}{\mathop{\mathrm{codim}}}
\newcommand{\rk}{\mathop{\mathrm{rk}}}
\newcommand{\X}{\Lambda}
\newcommand{\mult}{\mathop{\mathrm{mult}}\nolimits}
\newcommand{\res}{\mathop{\mathrm{res}}\nolimits}
\newcommand{\length}{\mathop{\mathrm{length}}}
\newcommand{\ch}{\mathop{\mathrm{ch}}\nolimits}
\newcommand{\Ll}{\mathcal{L}}
\newcommand{\Oo}{\mathcal{O}}
\newcommand{\eps}{\varepsilon}
\renewcommand{\phi}{\varphi}
\newcommand{\Ru}[1]{#1_{\text{u}}}
\newcommand{\SL}{\mathop{\mathrm{SL}}}
\newcommand{\divr}{\mathop{\mathrm{div}}}

\newtheorem{theorem}{Theorem}
\newtheorem{lemma}{Lemma}
\newtheorem*{Corollary}{Corollary}
\theoremstyle{definition}
\newtheorem{remark}{Remark}
\newtheorem*{Question}{Question}
\newtheorem{example}{Example}

\begin{document}

\begin{abstract}
The complexity of a homogeneous space $G/H$ under a reductive
group $G$ is by definition the codimension of general orbits in
$G/H$ of a Borel subgroup $B\subseteq G$. We give a
representation-theoretic interpretation of this number as the
exponent of growth for multiplicities of simple $G$-modules in
the spaces of sections of homogeneous line bundles on $G/H$. For this, we
show that these multiplicities are bounded from above by the
dimensions of certain Demazure modules. This estimate for
multiplicities is uniform, i.e., it depends not on~$G/H$, but
only on its complexity.
\end{abstract}

\maketitle

\section{Introduction}\label{intro}

Let $G$ be a connected reductive group over an algebraically
closed field $\kk$ of characteristic~0, and $H\subseteq G$ a
closed subgroup. Consider the homogeneous space $G/H$. Choose a
Borel subgroup $B\subseteq G$. By lower semicontinuity, general
$B$-orbits in $G/H$ have maximal dimension. The minimal
(=typical) codimension $c=c(G/H)$ of $B$-orbits is called the
\emph{complexity} of~$G/H$. (Clearly, it does not depend on the
choice of~$B$ since all Borel subgroups are conjugate.) By the
Rosenlicht theorem \cite[2.3]{inv.th},
$c(G/H)$ equals the transcendence degree of $\kk(G/H)^B$
over~$\kk$.

This numerical invariant plays an important r{\^o}le in the
geometry of~$G/H$. For instance, the class of homogeneous spaces
of complexity zero, called \emph{spherical} spaces, is
particularly nice \cite{LV(sph)}, \cite{var.sph}, \cite{comm}.
It includes many classical spaces, all symmetric spaces, etc.
Also the equivariant embedding theory of $G/H$ depends crucially
on its complexity, see \cite{LV}, \cite{LVT}.

In this note, we describe $c(G/H)$ in terms of representation
theory related to~$G/H$.

Let $\X=\X(B)$ be the weight lattice of~$B$, and
$\X_{+}\subseteq\X$ be the set of dominant weights. By
$V(\lambda)$ denote the simple $G$-module of
highest weight~$\lambda\in\X_{+}$. For any rational
$G$-module~$M$, let $\mult_{\lambda}M$ denote the multiplicity
of $V(\lambda)$ in~$M$.

It turns out that $c(G/H)$ characterizes the growth of
multiplicities in spaces of global sections of $G$-line bundles
over~$G/H$. Here is our main result:
\begin{theorem}\label{mult}
The complexity $c(G/H)$ is the minimal integer $c$ such that
$\mult_{\lambda}H^0(G/H,\Ll)=O(|\lambda|^c)$ over all
$\lambda\in\X_{+}$ and all $G$-line bundles $\Ll\to G/H$, where
$|\cdot|$ is any fixed norm on the vector space spanned by~$\X$.
Moreover, this estimate for multiplicities is uniform over all
$H\subseteq G$ such that $c(G/H)=c$. If $G/H$ is quasiaffine,
then it suffices to consider only $\mult_{\lambda}\kk[G/H]$.
\end{theorem}

A weaker version of Theorem~\ref{mult} under some restrictive
conditions (multiplicities in $\kk[G/H]$ for quasiaffine $G/H$
provided that $\kk[G/H]$ is finitely generated; no uniform
estimate) appeared in~\cite{mult}. The relation between
complexity and growth of multiplicities is known for quite a
time, see partial results in \cite[1.1]{t-dec}, \cite[2.4]{c&r},
\cite[1.3]{var.sph}.

We prove Theorem~\ref{mult} in Section~\ref{estimate}. The idea
is to embed the ``space of multiplicity'' in the dual of a
certain Demazure submodule in $V(\lambda)$, associated with an
element of length~$c$ in the Weyl group (Lemma~\ref{Demazure}).

In Section~\ref{c<2}, we justify the term ``complexity'' by
providing a much more precise information on multiplicities on
homogeneous spaces of complexity $\le1$. Actually, the spherical
case is well known \cite{mult(sph)} and is included in the text
only for convenience of the reader. In the case of complexity~1,
a formula similar to ours for multiplicities in $\kk[G/H]$
provided that it is a finitely generated unique factorization
domain, and $G/H$ is quasiaffine, was obtained
in~\cite[1.2]{t-dec}, see also \cite[2.4.19]{c&r}.

\subsection*{Acknowledgements}
This note was written during my stay at Institut Fourier in
spring 2003. I would like to thank this institution for
hospitality, and M.~Brion for invitation and for stimulating
discussions. Thanks are also due to I.~V.~Arzhantsev for some
helpful remarks.

\subsection*{Notation}
\begin{itemize}
\item The character lattice of an algebraic group $H$ is denoted
by $\X(H)$ and is written additively.
\item By $M^H$ we denote the set of $H$-fixed points in the set
$M$ acted on by a group~$H$. If $M$ is a vector space, and $H$
acts linearly, then $M_{\chi}=M^{(H)}_{\chi}$ is the
$H$-eigenspace of eigenweight $\chi\in\X(H)$.
\item Throughout the paper, $G$~is a connected reductive group,
and $B\subseteq G$ a fixed Borel subgroup. The semigroup
$\X_{+}$ of dominant weights is considered relative to~$B$. We
fix an opposite Borel subgroup $B^{-}$, a maximal torus $T=B\cap
B^{-}$, and consider the Weyl group $W$ of $G$ relative to~$T$.
\item By $\lambda^{*}$ denote the highest weight of the dual
$G$-module to the simple $G$-module $V(\lambda)$ of highest
weight $\lambda\in\X_{+}$.
\end{itemize}

\section{Upper bound for multiplicities}\label{estimate}

We begin with basic facts about line bundles on homogeneous
spaces. Every line bundle $\Ll\to G/H$ admits a
$G$-linearization, i.e., a fiberwise linear $G$-action
compatible with the projection onto the base, if we possibly
replace $G$ by a finite cover \cite[\S2]{loc.qpro}.
(Alternatively, one may replace $\Ll$ by its sufficiently big
tensor power.) Every $G$-line bundle is isomorphic to a
homogeneous bundle
$\Ll(\chi)=\Ll_{G/H}(\chi)=G\times^{H}\kk_{\chi}$, where $H$
acts on the fiber $\kk_{\chi}\simeq\kk$ via a character
$\chi\in\X(H)$. The bundle $\Ll(\chi)$ is trivial (regardless
the $G$-linearization) iff $\chi$ is the restriction of a
character of~$G$.

The space of global sections
$H^0(G/H,\Ll(\chi))\simeq\kk[G]^{(H)}_{-\chi}$ is a rational
$G$-module, where $H$ acts on $G$ by right translations, and $G$
by left translations. By the Frobenius reciprocity
\cite[I.3.3--3.4]{reps(AG)}, we have
$\mult_{\lambda}H^0(G/H,\Ll(\chi))=\dim
V(\lambda^{*})^{(H)}_{-\chi}$.

In particular, $H^0(G/B,\Ll(-\lambda))=V(\lambda^{*})$ whenever
$\lambda\in\X_{+}$, and~0, otherwise (the Borel--Weil theorem).

Observe that for any rational $G$-module $M$ we have
$\mult_{\lambda}M=\dim M^{(B)}_{\lambda}$,
$\forall\lambda\in\X_{+}$. In particular,
\begin{equation*}
\mult_{\lambda}H^0(G/H,\Ll(\chi))=
\dim\kk[G]^{(B\times H)}_{(\lambda,-\chi)}
\end{equation*}
The nonzero spaces $\kk[G]^{(B\times H)}_{(\lambda,-\chi)}$
represent complete linear systems of $B$-stable divisors
on~$G/H$, i.e., linear systems of pairwise rationally
equivalent $B$-stable divisors which cannot be enlarged by
adding new $B$-stable effective divisors. The respective
biweights $(\lambda,\chi)$ form a subsemigroup
$\Sigma=\Sigma(G/H)\subseteq\X(B\times H)$.  Two biweights
$(\lambda,\chi),(\lambda',\chi')\in\Sigma$ determine the same
linear system on $G/H$ iff $(\lambda,\chi)$ differs from
$(\lambda',\chi')$ by a twist of $G$-linearization, i.e., by
$(\eps|_B,-\eps|_H)$, $\eps\in\X(G)$.

We need a useful result, essentially due to Brion.
Replacing $H$ by a conjugate, we may assume that $\dim
B(eH)$ is maximal among all $B$-orbits, i.e., $\codim
B(eH)=c=c(G/H)$.
\begin{lemma}[{cf.~\cite[2.1]{param}}]\label{Schubert}
There exists a sequence of minimal parabolics
$P_1,\dots,P_c\supset B$ such that $\overline{P_c\cdots
P_1(eH)}=G/H$. The decomposition $w=s_1\cdots s_c$, where
$s_i\in W$ are the simple reflections corresponding to $P_i$, is
reduced, and $D_w=\overline{BwB}=P_1\cdots P_c$ is a
``Schubert subvariety'' in $G$ of dimension $c+\dim B$.
\end{lemma}
\begin{proof}
If $c>0$, then $\overline{B(eH)}$ is not $G$-stable, whence it
is not stabilized by some minimal parabolic $P_1\supset B$.
Since $P_1/B\simeq\PP^1$, the natural map $P_1\times^BB(eH)\to
G/H$ is generically finite, and $\codim P_1(eH)=c-1$.
Continuing in the same way, we construct a sequence of minimal
parabolics $P_1,\dots,P_c\supset B$ such that
$\overline{P_c\cdots P_1(eH)}=G/H$, i.e., $P_c\cdots P_1H$ is
dense in~$G$. The map $P_c\cdots P_1\times^BB(eH)\to G/H$ is
generically finite, hence $\dim P_c\cdots P_1=c+\dim B$, which
yields all the remaining assertions.
\end{proof}
\begin{remark}
The ``Schubert subvarieties'' $D_w$, $w\in W$, form a monoid w.r.t.\
the multiplication of sets in~$G$, called the
\emph{Richardson--Springer monoid}. It is generated by the $D_s$,
$s\in W$ a simple reflection, and is naturally identified with~$W$,
the defining relations being $s^2=s$ and the braid relations for~$W$.
The action of the Richardson--Springer monoid on the set of
$B$-stable subvarieties (appearing implicitly in
Lemma~\ref{Schubert}) is studied in~\cite{B-orb}.
\end{remark}

Let $v_{\lambda}\in V(\lambda)$ be a highest weight vector. The
$B$-submodule $V_w(\lambda)\subseteq V(\lambda)$ generated by
$wv_{\lambda}$, $w\in W$, is called a \emph{Demazure module}.
In the notation of Lemma~\ref{Schubert}, we have
$V_w(\lambda)=\langle D_wv_{\lambda}\rangle\simeq
H^0(S_w,\Ll_{G/B}(-\lambda))^{*}$, where $S_w=D_w/B$ is a
Schubert subvariety in $G/B$ of dimension~$c$. Indeed, the
restriction map $V(\lambda^{*})=H^0(G/B,\Ll(-\lambda))\to
H^0(S_w,\Ll(-\lambda))$ is surjective
\cite[II.14.15,\,e)]{reps(AG)}, and $D_wv_{\lambda}$ is the
affine cone over the image of $S_w$ under the map
$G/B\to\PP(V(\lambda))$.

\begin{lemma}\label{Demazure}
In the notation of Lemma~\ref{Schubert},
$\mult_{\lambda}H^0(G/H,\Ll(\chi))\le\dim V_w(\lambda)$.
\end{lemma}
\begin{proof}
We show that the pairing between $V(\lambda^{*})$ and
$V(\lambda)$ provides an embedding
$V(\lambda^{*})^{(H)}_{-\chi}\hookrightarrow V_w(\lambda)^{*}$.
Otherwise, if $v^{*}\in V(\lambda^{*})^{(H)}_{-\chi}$ vanishes on
$V_w(\lambda)$, then it vanishes on $D_wv_{\lambda}$,
i.e., $\langle D_{w}^{-1}v^{*}\rangle=\langle
Gv^{*}\rangle=V(\lambda^{*})$ vanishes at~$v_{\lambda}$, a
contradiction.
\end{proof}
\begin{remark}
A similar idea was used in \cite[2.4.18]{c&r} to obtain an upper
bound for multiplicities in coordinate algebras of homogeneous
spaces of complexity~1.
\end{remark}
\begin{remark}
The assertion of Lemma~\ref{Demazure} can be refined and viewed
in a more geometric context as follows. Consider the natural
$B$-equivariant proper map $\phi:D_w^{-1}/B\cap H\simeq
D_w^{-1}\times^BB(eH)\to G/H$, which is generically finite by
construction. For $\Ll=\Ll_{G/H}(\chi)$ we have
$\phi^{*}\Ll=\Ll_{D_w^{-1}/B\cap H}(\chi|_{B\cap H})$, and
$\phi^{*}:H^0(G/H,\Ll)\hookrightarrow
H^0(D_w^{-1}/B\cap H,\phi^{*}\Ll)$ gives rise to
\begin{multline*}
H^0(G/H,\Ll)^{(B)}_{\lambda}\hookrightarrow
H^0(D_w^{-1}/B\cap H,\phi^{*}\Ll)^{(B)}_{\lambda}=
\kk[D_w^{-1}]^{(B\times B\cap H)}_{(\lambda,-\chi)}\\
\simeq\kk[D_w]^{(B\cap H\times B)}_{(-\chi,\lambda)}=
H^0(S_w,\Ll(-\lambda))^{(B\cap H)}_{-\chi}=
{V_w(\lambda)^{*}}^{(B\cap H)}_{-\chi}
\end{multline*}
I am indebted to M.~Brion for this remark.
\end{remark}

Lemma~\ref{Demazure} applies to obtaining upper bounds for
multiplicities in branching to reductive subgroups,
cf.~\cite[Thm.\,2]{mult}.
\begin{Corollary}
If $L\subseteq G$ is a connected reductive subgroup, then
\begin{equation*}
\mult_{\mu}\res^G_LV(\lambda)\le\dim V_w(\lambda^{*})
\end{equation*}
for any two dominant weights $\lambda,\mu$ of $G,L$,
respectively, where $w\in W$ is provided by Lemma~\ref{Schubert}
for $H$ equal to a Borel subgroup of~$L$. Similarly,
\begin{equation*}
\length\res^G_LV(\lambda)\le\dim V_w(\lambda^{*})
\end{equation*}
where $w\in W$ corresponds to $H$ equal to a maximal unipotent
subgroup of~$L$. (Here $\length$ is the number of simple factors
in an $L$-module.)
\end{Corollary}
\begin{proof}
Just note that
\begin{equation*}
\mult_{\mu}\res^G_LV(\lambda)=\dim
V(\lambda)^{(H)}_{\mu}=\mult_{\lambda^{*}}H^0(G/H,\Ll(-\mu))
\end{equation*}
in
the first case, and $\length\res^G_LV(\lambda)=\dim
V(\lambda)^H=\mult_{\lambda^{*}}\kk[G/H]$ in the second case,
and then apply Lemma~\ref{Demazure}.
\end{proof}

\begin{example}\label{uni(Levi)}
Let $P\supseteq B^{-}$ be a parabolic subgroup with the Levi
decomposition $P=L\rightthreetimes\Ru{P}$, $L\supseteq T$. By
$w_L$ denote the longest element in the Weyl group of~$L$, and
consider the decomposition $w_G=w_Lw^L$. Let $H$ be the
unipotent radical of $B^{-}\cap L$. Then we may take $w=w^L$ and
obtain $\length\res^G_LV(\lambda)\le\dim V_{w^L}(\lambda^{*})$.
\end{example}

\begin{proof}[Proof of Theorem~\ref{mult}]
Recall the character formula for Demazure modules
\cite[II.14.18,\,b)]{reps(AG)}:
\begin{equation*}
\ch_TV_w(\lambda)=\frac{1-e^{-\alpha_1}s_1}{1-e^{-\alpha_1}}\dots
\frac{1-e^{-\alpha_c}s_c}{1-e^{-\alpha_c}}\ e^{\lambda}
\end{equation*}
where $e^{\mu}$ is the monomial in the group algebra $\ZZ[\X]$
corresponding to $\mu\in\X$, $\alpha_i$~are the simple roots
defining~$P_i$, and $s_i\in W$ are the respective simple
reflections acting on $\ZZ[\X]$ in a natural way. One easily
computes
\begin{equation*}
\frac{1-e^{-\alpha_i}s_i}{1-e^{-\alpha_i}}\ e^{\mu}=
\begin{cases}
e^{\mu}(1+e^{-\alpha_i}+\dots+
e^{-\langle\mu,\alpha_i^{\vee}\rangle\alpha_i}), &
\langle\mu,\alpha_i^{\vee}\rangle\ge0, \\
0, & \langle\mu,\alpha_i^{\vee}\rangle=-1, \\
-e^{\mu}(e^{\alpha_i}+e^{2\alpha_i}+\dots+
e^{|1+\langle\mu,\alpha_i^{\vee}\rangle|\alpha_i}), &
\langle\mu,\alpha_i^{\vee}\rangle\le-2,
\end{cases}
\end{equation*}
$\forall i=1,\dots,c,\ \forall\mu\in\X$, where $\alpha_i^{\vee}$
is the respective simple coroot. It is then easy to deduce that
$\dim V_w(\lambda)=O(|\lambda|^c)$, and Lemma~\ref{Demazure}
yields the desired estimate in Theorem~\ref{mult}.

Alternatively, observing that $S_w$ is a projective variety of
dimension~$c$, one may deduce that $\dim H^0(S_w,\Ll(-\lambda))$
grows no faster than~$|\lambda|^c$ as follows.

Without loss of generality we may assume $G$ to be semisimple
and simply connected. Let $\omega_1,\dots,\omega_l$ be the
fundamental weights of~$G$. Put
$X=\PP_{S_w}(\Ll(\omega_1)\oplus\dots\oplus\Ll(\omega_l))$, a
projective space bundle over $S_w$ with fiber~$\PP^{l-1}$. Let
$\Oo_X(1)$ be the antitautological line bundle over $X$, and
$\pi:X\to S_w$ the projection map. Then
\begin{align*}
\pi_{*}\Oo_X(k)&=\bigoplus_{k_1+\dots+k_l=k,\ k_j\ge0}
\Ll(-k_1\omega_1-\dots-k_l\omega_l) \\
R^i\pi_{*}\Oo_X(k)&=0,\qquad \forall i>0
\end{align*}
Now $R_w=\bigoplus_{k\ge0}H^0(X,\Oo(k))=
\bigoplus_{\lambda\in\Lambda_{+}}H^0(S_w,\Ll(-\lambda))$ is a
$\Lambda$-graded algebra, and the quotient field of $R_w$ is a
monogenic transcendental extension of~$\kk(X)$, hence it has
transcendence degree $n=c+l$. Furthermore, $R_w$~is finitely
generated.  Indeed, $R_w$ is a quotient of the $G$-algebra
$R=\bigoplus_{\lambda\in\Lambda_{+}}H^0(G/B,\Ll(-\lambda))$.
The multihomogeneous components of~$R$ are the simple
$G$-modules~$V(\lambda^{*})$, hence $R$ is generated by its
components with $\lambda=\omega_1,\dots,\omega_l$. (The spectrum
of $R_w$ is the so called \emph{multicone} over the Schubert
variety $S_w$, studied in~\cite{multicones}. For instance, this
multicone has rational singularities.) Finally, a general
property of finitely generated multigraded algebras implies that
$\dim H^0(S_w,\Ll(-\lambda))=O(|\lambda|^{n-\rk\Lambda})=
O(|\lambda|^c)$.

This estimate is uniform in~$H$, because there are finitely many
choices for~$w$. It remains to show that the exponent $c$ cannot
be made smaller.

Let $f_1,\dots,f_c$ be a transcendence base of $\kk(G/H)^B$.
There exist a $G$-line bundle $\Ll$ and $B$-eigenvectors
$\sigma_0,\dots,\sigma_c\in H^0(G/H,\Ll)$ of the same weight
$\lambda$ such that $f_i=\sigma_i/\sigma_0$, $\forall i=1,\dots,c$.
(Indeed, $\Ll$~and $\sigma_0$ may be determined by a sufficiently big
$B$-stable effective divisor majorizing the poles of all~$f_i$.)

Consider the graded algebra $R=\bigoplus_{n\ge0}R_n$,
$R_n=H^0(G/H,\Ll^{\otimes n})^{(B)}_{n\lambda}$. Clearly,
$\sigma_0,\dots,\sigma_c$ are algebraically independent in~$R$,
whence
\begin{equation*}
\mult_{n\lambda}H^0(G/H,\Ll^{\otimes n})=\dim
R_n\ge\binom{n+c}{c}\sim n^c
\end{equation*}
This proves our claim.

Finally, if $G/H$ is quasiaffine, then there even exist
$\sigma_0,\dots,\sigma_c\in\kk[G/H]$ with the same properties.
\end{proof}

\begin{remark}
For a given~$H$, the estimate of the multiplicity by $\dim
V_w(\lambda)$ may be not sharp. However, it is natural to ask
whether it is a sharp uniform estimate over all homogeneous
spaces with given complexity. More precisely, we formulate the
following
\begin{Question}
Given an element $w\in W$ of length~$c$, does there exist a
subgroup $H\subseteq G$ such
that $\mult_{\lambda}H^0(G/H,\Ll(\chi))=\dim V_w(\lambda)$ for
sufficiently general $(\lambda,\chi)\in\Sigma(G/H)$?
\end{Question}
\end{remark}

\begin{example}\label{uni.rad}
In the notation of Example~\ref{uni(Levi)}, put $H=\Ru{P}$. Then
always $\chi=0$, so that $H^0(G/H,\Ll(\chi))=\kk[G/H]$. We may
take $w=w_L$. Then $V(\lambda^{*})^H$ is a simple $L$-module of
lowest weight~$-\lambda$, and $V_w(\lambda)$ is the dual
$L$-module of highest weight~$\lambda$. It follows that
$\mult_{\lambda}\kk[G/H]=\dim V_w(\lambda)$.
\end{example}

\section{Case of small complexity}\label{c<2}

Homogeneous spaces of complexity $\le1$ are distinguished among
all homogeneous spaces by their nice behaviour. For instance,
they have a well-developed equivariant embedding theory
\cite[8--9]{LV}, \cite[2--5]{LVT}. There are also more explicit
formul{\ae} for multiplicities in this case.

\begin{theorem}\label{mult(c<2)}
In the above notation,
\begin{enumerate}
\item\label{mult(c=0)} If $c(G/H)=0$, then
$\mult_{\lambda}H^0(G/H,\Ll(\chi))=1$,
$\forall(\lambda,\chi)\in\Sigma$.
\item\label{mult(c=1)} If $c(G/H)=1$, then there exists a pair
$(\lambda_0,\chi_0)\in\Sigma$, unique up to a shift by
$(\eps|_B,-\eps|_H)$, $\eps\in\X(G)$, such that
\begin{equation*}
\mult_{\lambda}H^0(G/H,\Ll(\chi))=n+1
\end{equation*}
where $n$ is the maximal integer such that
$(\lambda,\chi)-n(\lambda_0,\chi_0)\in\Sigma(G/H)$.
\end{enumerate}
\end{theorem}
\begin{proof}
The assertion is well known in the case $c=0$, and we prove it
just to keep the exposition self-contained. Assuming the
contrary yields two non-proportional $B$-eigenvectors
$\sigma_0,\sigma_1\in H^0(G/H,\Ll(\chi))$ of the same
weight~$\lambda$. Hence $f=\sigma_1/\sigma_0\in\kk(G/H)^B$,
$f\ne\text{const}$, a contradiction.

In the case $c=1$, we have $\kk(G/H)^B\simeq\kk(\PP^1)$ by the
L{\"u}roth theorem. Consider the respective rational map
$\pi:G/H\dasharrow\PP^1$, whose general fibers are (the closures
of) general $B$-orbits. By a standard argument, $\pi$~is given by
two $B$-eigenvectors $\sigma_0,\sigma_1\in H^0(G/H,\Ll(\chi_0))$
of the same weight $\lambda_0$ for a certain
$(\lambda_0,\chi_0)\in\Sigma$. Moreover, $\sigma_0,\sigma_1$
are algebraically independent, and each $f\in\kk(G/H)^B$ can be
represented as a homogeneous rational fraction in
$\sigma_0,\sigma_1$ of degree~0.

Now put $(\mu,\tau)=(\lambda,\chi)-n(\lambda_0,\chi_0)$, fix
$\sigma_{\mu}\in H^0(G/H,\Ll(\tau))^{(B)}_{\mu}$, and take any
$\sigma_{\lambda}\in H^0(G/H,\Ll(\chi))^{(B)}_{\lambda}$. Then
$f=\sigma_{\lambda}/\sigma_0^n\sigma_{\mu}\in\kk(G/H)^B$, whence
$f=F_1/F_0$ for some $m$-forms $F_0,F_1$ in $\sigma_0,\sigma_1$.
We may assume the fraction to be reduced and decompose
$F_1=L_1\cdots L_m$, $F_0=M_1\cdots M_m$, as products of linear
forms, with all $L_i$ distinct from all~$M_j$. Then
$\sigma_{\lambda}M_1\cdots M_m=\sigma_{\mu}\sigma_0^nL_1\cdots
L_m$.

Being fibers of~$\pi$, the divisors of $\sigma_0$, $L_i$, $M_j$
on $G/H$ either coincide or have no common components. By the
maximality of~$n$, the divisor of $\sigma_{\mu}$ does not
majorize any one of $M_j$. Therefore $M_1=\dots=M_m=\sigma_0$,
$m\le n$, and $\sigma_{\lambda}/\sigma_{\mu}$ is an $n$-form in
$\sigma_0,\sigma_1$. The assertion follows.
\end{proof}

\begin{example}
Let $G=\SL_3$, $H=T$ (the diagonal torus), $B$ be the
upper-triangular subgroup. The space $G/H$ can be regarded as
the space of ordered triangles in~$\PP^2$, i.e.,
the set of all ordered triples
$p=(p_1,p_2,p_3)\in\PP^2\times\PP^2\times\PP^2$ such that
$p_1,p_2,p_3$ do not lie on one line.
Let $\ell_i\subset\PP^2$ denote the line joining $p_j$
and~$p_k$, where $(i,j,k)$ is a cyclic permutation of $(1,2,3)$.
By $p_0$ denote the $B$-fixed point in $\PP^2$, and by $\ell_0$
the $B$-stable line.

There are the following $B$-stable prime divisors on~$G/H$:
\begin{xalignat*}{3}
D_i&=\{p\mid p_i\in\ell_0\}=\divr g_{3i},&
\lambda_i&=\omega_2,\ \chi_i=-\eps_i \\
D_i'&=\{p\mid p_0\in\ell_i\}=\divr\Delta_i,\quad
\Delta_i=\begin{vmatrix}
g_{2j} & g_{2k} \\
g_{3j} & g_{3k}
\end{vmatrix},&
\lambda'_i&=\omega_1,\ \chi'_i=\eps_i \\
D_t&=\overline{B\cdot p(t)}=
\divr(g_{32}\Delta_2+tg_{33}\Delta_3),&
\lambda_0&=\omega_1+\omega_2,\ \chi_0=0
\end{xalignat*}
where $i=1,2,3$, $t\in\PP^1\setminus\{0,1,\infty\}$, and the
vertices of the triangle $p(t)$ are: $p_1(t)=(0:0:1)$,
$p_2(t)=(0:1:1)$, $p_3(t)=(1:t:1)$. Here $g_{ij}$ are matrix
entries of $g\in G$, and $H$-semiinvariant polynomials in
$g_{ij}$ are regarded as sections of $G$-line bundles on~$G/H$.
We also indicate their biweights $(\lambda,\chi)\in\Sigma$,
denoting by $\omega_i$ the fundamental weights, and by $\eps_i$
the diagonal matrix entries of~$H$. Observe that
$g_{31}\Delta_1+g_{32}\Delta_2+g_{33}\Delta_3=0$.

It follows that $c(\SL_3/T)=1$. Now it is an easy combinatorial
exercise to deduce from Theorem~\ref{mult(c<2)}(\ref{mult(c=1)})
that $\mult_{\lambda}H^0(\SL_3/T,\Ll(\chi))=n+1$, where
\begin{equation*}
n=\frac{k_1+k_2}2-\frac16\sum_{i=1}^3
|k_1-k_2+2l_i-l_j-l_k|
\end{equation*}
whenever $(\lambda,\chi)\in\Sigma$,
$\lambda=k_1\omega_1+k_2\omega_2$,
$\chi=l_1\eps_1+l_2\eps_2+l_3\eps_3$; and
$(\lambda,\chi)\in\Sigma$ whenever $k_1-k_2\equiv l_1+l_2+l_3
\pmod{3}$ and $n\ge0$.
\end{example}

\end{document}